\documentclass[11 pt,a4paper,twoside,reqno]{amsart} 
\usepackage{amsfonts,amssymb,amscd,amsmath,enumerate,verbatim,calc} 
 
\usepackage{setspace}
\usepackage{amsthm,mathrsfs,textcomp}
\usepackage[utf8]{inputenc}	
\usepackage[T1]{fontenc}	
\usepackage{mathtools}
\usepackage[super]{nth}
\usepackage[all,arc]{xy}
\usepackage{hyperref}
\usepackage{commath}
\usepackage[english]{babel}
\newcommand\numberthis{\addtocounter{equation}{1}\tag{\theequation}}

\theoremstyle{plain}

\theoremstyle{definition}

\numberwithin{equation}{section}

\newcommand {\Mathstu}{\hspace*{-0.6 cm}The Mathematics Student\hfill \hspace{2 in}\hfill ISSN: 0025-5742\\ Vol. 93, Nos. 1-2,  January-June (2024), 24--29}
\begin{document}
	\Mathstu
	\vspace{0.7 cm}
	
	\title
	[Proofs of Two Formulas of Vladeta Jovovic]
	{\scshape \large \bf {Proofs of Two Formulas of Vladeta Jovovic}}
	\vspace*{-0.4 cm}
	\author [Aritram Dhar] {Aritram Dhar}
	\thanks{\raggedright{\hspace*{-.21 in}
\hspace{0.24cm}			{2010 Mathematics Subject Classification}: 05A17, 05A19, 11P81, 11P82    }  \\ [0.4 cm]
\hspace{5.5cm}		\raggedleft \textit{\large \copyright\,\,{Indian Mathematical Society, 2024}}}
	
	\maketitle
	\vspace*{-0.9 cm}

\begin{abstract}
In this paper, we first provide an analytic and a bijective proof of a formula stated by Vladeta Jovovic in the OEIS sequence A117989. We also provide a bijective proof of another interesting result stated by him on the same page concerning integer partitions with fixed differences between the largest and smallest parts.
\end{abstract}

\maketitle

\section{Introduction}\label{d1}
A partition $\pi$ of a positive integer $n$ is a non-increasing sequence of natural numbers $\lambda_1 \ge \lambda_2 \ge \ldots \ge \lambda_r$ such that $\displaystyle{\sum_{i=1}^{r}\lambda_i = n}$. $\lambda_1, \lambda_2, \ldots, \lambda_r$ are called the \textit{parts} of the partition $\pi$. We call $\lambda_1$ and $\lambda_r$ to be the \textit{largest} and \textit{smallest} parts of the partition $\pi$ respectively and denote $p(n)$ to be the number of (unrestricted) partitions of $n$. For example, $p(5) = 7$ where the seven partitions of $5$ are $5$, $4+1$, $3+2$, $3+1+1$, $2+2+1$, $2+1+1+1$, and $1+1+1+1+1$. Their corresponding largest parts are $5$, $4$, $3$, $3$, $2$, $2$, and $1$ respectively and smallest parts are $5$, $1$, $2$, $1$, $1$, $1$, and $1$ respectively.\par Throughout the paper, we consider $\abs{q}<1$ and adopt the usual notation for the conventional $q$-Pochammer symbols:
$$(a)_n = (a;q)_n := \prod\limits_{k=0}^{n-1}(1-aq^k),$$
$$(a)_{\infty} = (a;q)_{\infty} := \lim_{n\to\infty}(a;q)_n.$$
We define $(a)_n$ for all real numbers $n$ by
$$(a)_n := \dfrac{(a)_{\infty}}{(aq^n)_{\infty}}.$$
\par By convention, we take $(a)_0 = 1$ and $(0)_{\infty} = 1$.

\section{Main Formulas}  
Our starting point is the sequence A$117989$ in the On-Line Encyclopedia of Integer Sequences \cite{3}. The sequence in question, $a(n)$, counts the number of partitions of $n$ where the smallest part occurs at least twice. For example, $a(6) = 7$ where the relevant partitions are $4+1+1$, $3+3$, $3+1+1+1$, $2+2+2$, $2+2+1+1$, $2+1+1+1+1$, and $1+1+1+1+1+1$.\par Also, on the page of A$117989$, we find the following formula of $a(n)$ given by Vladeta Jovovic:\newline \quad \newline \textbf{Formula $1$:} (Vladeta Jovovic, July 21 2006) \begin{equation} \label{eq1}
    a(n) = 2p(n) - p(n+1) \,\,\, \forall \,\, n\ge 1
\end{equation}\textit{where $p(n)$ denotes the partition function}.\newline \quad \par We denote $b(n) = 2p(n)-p(n+1)$. By (\ref{eq1}), we have $a(6) = 7 = 22-15 = 2p(6)-p(7) = b(6)$.\par We also find another interesting result on the same page posted by Vladeta Jovovic which states:\newline \quad \newline \textbf{Formula $2$:} (Vladeta Jovovic, May 09 2008) \begin{equation} \label{eq2}
    a(n) = p(2n,n) \,\,\, \forall \,\, n\ge 1
\end{equation}\textit{where $p(2n,n)$ denotes the number of partitions of $2n$ with fixed difference equal to $n$ between the largest and smallest parts}.\newline \quad \par In sections \ref{s3} and \ref{s4}, we give a q-theoretic and a bijective proof of (\ref{eq1}) respectively and finally we provide a bijective proof of (\ref{eq2}) in section \ref{s5}.

\section{Analytic Proof of Formula 1}\label{s3}
In this section, we prove (\ref{eq1}) using generating functions and elementary infinite series-product identities from the theory of $q$-hypergeometric series.\par To begin with, we define the generating functions for $a(n)$ and $b(n)$ to be\begin{equation}
    A(q) := \displaystyle{\sum_{n=1}^{\infty}a(n)q^n}\quad \text{and}\quad 
    B(q) := \displaystyle{\sum_{n=1}^{\infty}b(n)q^n}
\end{equation}respectively. We then have\begin{align*}
    A(q)&=\displaystyle{\sum_{k=1}^{\infty}q^{k+k}(1+q^k+q^{2k}+\cdots)(1+q^{k+1}+q^{2(k+1)}+\cdots)\cdots}
    \\
    &=\displaystyle{\sum_{k=1}^{\infty}\frac{q^{2k}}{(q^k)_{\infty}}}
    \\
    &=\frac{1}{(q)_{\infty}}\displaystyle{\sum_{k=1}^{\infty}(q)_{k-1}q^{2k}}
    \\
    &=\frac{q^2}{(q)_{\infty}}\displaystyle{\sum_{k=0}^{\infty}\frac{(q)_k(q)_kq^{2k}}{(q)_k}}
    \\
    &=\frac{q^2(q^3)_{\infty}}{(q^2)_{\infty}}\displaystyle{\sum_{k=0}^{\infty}\frac{(q^2)_kq^k}{(q)_k(q^3)_k}}\numberthis\label{eq3}
    \\
    &= \displaystyle{\sum_{k=0}^{\infty}\frac{q^{k+2}}{(q)_k(1-q^{k+2})}}
    \\
    &=\displaystyle{\sum_{k=0}^{\infty}\frac{q^{k+2}(1-q^{k+1})}{(q)_{k+2}}}
    \\
    &=\displaystyle{\sum_{k=0}^{\infty}\frac{q^{k+2}}{(q)_{k+2}}-\sum_{k=0}^{\infty}\frac{q^{2k+3}}{(q)_{k+2}}}
    \\
    &=\displaystyle{\sum_{k=0}^{\infty}\frac{q^{k+2}}{(q)_{k+2}}-\frac{1}{q}\sum_{k=0}^{\infty}\frac{(q^2)^{k+2}}{(q)_{k+2}}}
    \\
    &=\bigg(\frac{1}{(q)_{\infty}}-\frac{1}{1-q}\bigg)-\frac{1}{q}\bigg(\frac{1}{(q^2)_{\infty}}-\frac{1-q+q^2}{1-q}\bigg)\numberthis\label{eq4}
    \\
    &=\frac{1}{(q)_{\infty}}-\frac{1}{1-q}-\frac{1-q}{q(q)_{\infty}}+\frac{1-q+q^2}{q(1-q)}
    \\
    &=\frac{2}{(q)_{\infty}}-\frac{1}{q(q)_{\infty}}+\frac{1}{q}-1
\end{align*} which is equal to $B(q)$, the generating function for $b(n) = 2p(n)-p(n+1)$ $\forall$ $n\ge 1$.\par Note that (\ref{eq3}) follows by replacing $a = q$, $b = q$, $c = 0$, and $t = q^2$ in Heine's transformation \cite{1}, p.$19$, Corollary $2.3$] and (\ref{eq4}) follows by replacing $a=0$, $t=q$ and $a=0$, $t=q^2$ respectively in Cauchy's identity [\cite{1}, p.$17$, Theorem $2.1$].\qed

\section{Bijective Proof of Formula 1}\label{s4}
In this section, we provide a bijective proof of (\ref{eq1}).\par From (\ref{eq1}), we have\begin{align*}
    a(n) &= 2p(n) - p(n+1)\\ &= p(n) - (p(n+1) - p(n))
\end{align*} which implies\begin{equation}\label{eq5}
    p(n) - a(n) = p(n+1) - p(n).
\end{equation}\par Let us now define $c(n) = p(n) - a(n)$ and $d(n) = p(n+1) - p(n)$. Thus, from (\ref{eq5}), it suffices to prove that $c(n) = d(n)$ $\forall$ $n\ge 1$.\par We note that $c(n)$ denotes the number of partitions of $n$ where the smallest part occurs exactly once. This follows straightforward from the definition of $a(n)$. We also note that $d(n)$ denotes the number of partitions of $n+1$ which do not contain $1$ as a part because every partition of $n+1$ which contains $1$ as a part can be obtained by adjoining $1$ as a part to every partition of $n$.\par Let $\mathcal{C}_n$ be the set of all partitions of $n$ where the smallest part occurs exactly once and $\mathcal{D}_n$ be the set of all partitions of $n+1$ not containing $1$ as a part. So, $\#\mathcal{C}_n = c(n)$ and $\#\mathcal{D}_n = d(n)$. Thus, it is clear that we will now produce a bijection between the sets $\mathcal{C}_n$ and $\mathcal{D}_n$ to obtain the desired result.\par Firstly, we consider a partition $\pi\in \mathcal{C}_n$ and consider two cases pertaining to $\pi$: If $1$ is a part of $\pi$, since it is the smallest part, it occurs exactly once. Now, add $1$ to the $1$ already in $\pi$ to get a new partition $\pi^{\prime}\in \mathcal{D}_n$ whose smallest part now is $2$. Hence, $\pi^{\prime}$ does not contain $1$ as a part. Now, if $1$ is not a part of $\pi$, then the smallest part of $\pi$ is greater than or equal to $2$ and hence it does not contain $1$ as a part. On adding $1$ to the smallest part of $\pi$, we get a new partition $\pi^{\prime}$ of $n+1$ which does not contain $1$ as a part. Hence, $\pi^{\prime}\in \mathcal{D}_n$.\par Now, we consider a partition $\pi^{\prime}\in \mathcal{D}_n$. Thus, $\pi^{\prime}$ does not contain $1$ as a part which implies that the smallest part of $\pi^{\prime}$ is greater than or equal to $2$. Again, we consider two cases concerning $\pi^{\prime}$: If the smallest part of $\pi^{\prime}$ occurs exactly once, we subtract $1$ from it to get a new partition $\pi\in \mathcal{C}_n$ and we are done. On the other hand, if the smallest part of $\pi^{\prime}$ occurs at least twice, subtract $1$ from any one of the smallest parts to get a new partition $\pi\in \mathcal{C}_n$.\par Thus, the process is reversible and hence $\mathcal{C}_n$ is bijection with $\mathcal{D}_n$ $\forall$ $n\ge 1$. So, we have our desired result.\qed

\section{Bijective Proof of Formula 2}\label{s5}
In this section, we provide a bijective proof of (\ref{eq2}).\par Let $\mathcal{A}_n$ be the set of all partitions of $n$ where the smallest part occurs at least twice and $\mathcal{F}_n$ be the set of all partitions of $2n$ where the difference between the largest and smallest parts is equal to $n$. So, $\#\mathcal{A}_n = a(n)$ and $\#\mathcal{F}_n = p(2n,n)$. Now, we will provide a bijection between the sets $\mathcal{A}_n$ and $\mathcal{F}_n$ to show that $a(n) = p(2n,n)$.\par Firstly, we consider a partition $\pi\in \mathcal{A}_n$. Then, we add $n$ to any one of the smallest parts of $\pi$ (since the smallest part of $\pi$ occurs at least twice) to get a new partition $\pi^{\prime}$. $\pi^{\prime}\in \mathcal{F}_n$ because the largest part of $\pi^{\prime}$ now is equal to $n+$ the smallest part of $\pi$ and the smallest part of $\pi^{\prime}$ is equal to the smallest part of $\pi$ thus making the difference equal to $n$.\par For the other way, we now consider a partition $\pi^{\prime}\in \mathcal{F}_n$. Note that the largest part of $\pi^{\prime}$ occurs exactly once. We then subtract $n$ from the largest part of $\pi^{\prime}$ to get a new partition $\pi$ whose smallest part is equal to the smallest part of $\pi^{\prime}$ and consider two cases pertaining to $\pi^{\prime}$: If the smallest part of $\pi^{\prime}$ occurs at least twice, we are done, i.e., $\pi\in \mathcal{A}_n$ since subtracting $n$ from the largest part of $\pi^{\prime}$ does not affect the frequency of the smallest part of $\pi$ ($=$ the smallest part of $\pi^{\prime}$) which still remains at least $2$. Lastly, if the smallest part of $\pi^{\prime}$ occurs exactly once, subtracting $n$ from the largest part of $\pi^{\prime}$ makes it equal to the the smallest part of $\pi^{\prime}$ and thus, the new partition $\pi$ that we obtain has smallest part occuring at least twice since smallest parts of $\pi^{\prime}$ and $\pi$ are equal. Hence, $\pi\in \mathcal{A}_n$.\par Thus, the process is reversible and hence, we have a bijection between $\mathcal{A}_n$ and $\mathcal{F}_n$ giving our desired result.\qed   

\section{Conclusion}
Thus, we have $p(2n,n) = a(n) = 2p(n) - p(n+1)$ $\forall$ $n\ge 1$. We speculate that there is an interesting proof of (\ref{eq2}) using the generating function approach. Although in \cite{2}, Andrews, Beck, and Robbins gave the generating function for $p(n,t)$ (which is the number of partitions of $n$ with fixed difference equal to $t$ between the largest and smallest parts), the generating function of $p(2n,n)$ does not follow straightforward.\par In the same spirit, we define $G_m(q) := \displaystyle{\sum_{n=1}^{\infty}a_m(n)q^n}$ where $a_m(n)$ denotes the number of partitions of $n$ where the smallest part occurs at least $m$ times. On the page of A$117989$, we see that $G_m(q) = \displaystyle{\sum_{k=1}^{\infty}\frac{q^{mk}}{(q^k)_{\infty}}}$. It will be very interesting to see if there is a closed formula analogous to (\ref{eq1}) for $a_m(n)$ $\forall$ $m\ge 3$ and if there exists such a formula, then it would be nice to provide a combinatorial proof of it. 

\noindent \textbf{Acknowledgement:}
The author would like to thank George E. Andrews for suggesting him to prove the two formulas of Vladeta Jovovic which came up during an ongoing project with him.

\bigskip

\bigskip

\noindent \textsc{Aritram Dhar\\ Department of Mathematics\\ University of Florida, \\ Gainesville, FL 32611.}\\
E-mail: aritramdhar@ufl.edu

\end{document}